\documentclass[12pt,a4paper]{article}
\usepackage[utf8]{inputenc}
\usepackage[english]{babel}
\usepackage{amsmath,amsfonts,amssymb,amsthm}
\usepackage{tikz}
\usetikzlibrary{cd}
\title{Residual finiteness of extensions of
arithmetic subgroups of $\SU(d,1)$ with cusps}
\author{Richard M. Hill\\
	\small{University College London}\\
	\small{\texttt{r.m.hill@ucl.ac.uk}}}
\date{}

\usepackage{amsmath,amssymb,amsfonts,amsthm}
\usepackage{hyperref}

\usepackage[top=2.5cm, left=2cm, right=2cm, bottom=4.0cm]{geometry}

\newtheorem{lemma}{Lemma}
\newtheorem{proposition}{Proposition}
\newtheorem{theorem}{Theorem}

\newtheorem*{thm}{Theorem}
\newcommand{\optional}[1]{}

\DeclareMathOperator{\ab}{ab}
\DeclareMathOperator{\BS}{{BS}}
\DeclareMathOperator{\cts}{{cts}}
\DeclareMathOperator{\compact}{compact}

\DeclareMathOperator{\Cong}{Cong.}

\DeclareMathOperator{\Hom}{Hom}

\DeclareMathOperator{\meas}{meas}

\DeclareMathOperator{\stable}{stable}

\DeclareMathOperator{\Rest}{Rest}

\DeclareMathOperator{\SL}{SL}

\DeclareMathOperator{\SU}{SU}

\DeclareMathOperator{\U}{U}

\newcommand{\A}{\mathbb{A}}
\newcommand{\C}{\mathbb{C}}

\newcommand{\G}{\mathbb{G}}

\newcommand{\Proj}{\mathbb{P}}
\newcommand{\Q}{\mathbb{Q}}
\newcommand{\R}{\mathbb{R}}
\newcommand{\Z}{\mathbb{Z}}

\newcommand{\cH}{\mathcal{H}}
\newcommand{\cO}{\mathcal{O}}

\newcommand{\SUd}{\SU(d,1)}

\begin{document}
\maketitle

\begin{abstract}
	Let $\Gamma$ be an arithmetic subgroup of $\SU(d,1)$ with cusps, and let $X_\Gamma$ be the associated locally symmetric space.
	In this paper we investigate the pre-image of $\Gamma$ in the covering groups of $\SU(d,1)$.
	Let $H^\bullet_!(X_\Gamma,\C)$ be the inner cohomology, i.e. the image
	in $H^\bullet(X_\Gamma,\C)$ of the compactly supported cohomology.
	We prove that if the first inner cohomology group $H^1_!(X_\Gamma,\C)$ is non-zero
	then the pre-image of $\Gamma$ in each connected cover of $\SU(d,1)$ is residually finite.
	At the end of the paper we give an example of an arithmetic subgroup $\Gamma$ satisfying the criterion $H^1_!(X_\Gamma,\C) \ne 0$.



%\paragraph{Availability of data and material:}
%Data sharing not applicable to this article as no datasets were generated or analysed during the current study.

\paragraph{Competing Interests:}
None.

\end{abstract}

\section{Introduction}

Let $\Gamma$ be an arithmetic subgroup of $\SU(d,1)$ for some $d \ge 2$.
The universal cover $\widetilde{\SU(d,1)}$ of $\SU(d,1)$ is
an infinite cyclic cover, so that we have a central extension
\[
	1 \to \Z \to \widetilde{\SU(d,1)} \to \SU(d,1) \to 1.
\]
Furthermore, for every positive integer $n$, there is a unique connected $n$-fold cover of $\SU(d,1)$, which is isomorphic to $\widetilde{\SU(d,1)}/n\Z$.
Let $\tilde \Gamma^{(n)}$ be the pre-image of $\Gamma$ in the connected $n$-fold cover, and let $\tilde \Gamma$ be the pre-image of $\Gamma$ in the universal cover.
In this paper, we shall investigate whether the groups
$\tilde \Gamma$ and $\tilde\Gamma^{(n)}$ are residually finite.
There are two motivations for studying this question.
\begin{enumerate}
\item
It is a famous open question whether every word-hyperbolic group is residually finite. If all such groups are indeed residually finite, then every $\tilde\Gamma^{(n)}$ must be residually finite. This would imply that $\tilde\Gamma$ is also residually finite.
\item
If $\tilde\Gamma^{(n)}$ is residually finite, then for every sufficiently large integer $m$, there exist modular forms on $\SU(d,1)$ of weight $\frac{m}{n}$ whose level is a subgroup of finite index in $\Gamma$.
The existence of such modular forms is discussed in \cite{hill}, and some examples of the forms (of weight $\frac{1}{3}$) have recently been described in \cite{freitag hill}.
\end{enumerate}

The question of whether $\tilde \Gamma$ and $\tilde\Gamma^{(n)}$ are residually finite has recently been studied in \cite{stover toledo}, \cite{hill}
\cite{stover toledo 2} and \cite{freitag hill}.
In both \cite{hill}, and
\cite{stover toledo 2}, it is is shown (independently) that for a certain class of cocompact arithmetic subgroups $\Gamma$, the groups $\tilde\Gamma$ and $\tilde\Gamma^{(n)}$ are residually finite.
In the current paper,  we give some evidence that a similar result might be true when $\Gamma$ has cusps.

To put our result in context, we first recall a theorem from
\cite{hill} and \cite{stover toledo 2}.
Let $k$ be a CM field of degree $[k:\Q]=2e$.
Choose a $(d+1)\times (d+1)$ Hermitian matrix $J$ with entries in $k$, such that
\begin{itemize}
	\item
	The matrix $J$ has signature is $(d,1)$ at one of the complex places of $k$;
	\item
	For each of the other $e-1$ complex places of $k$,
	the matrix $J$ is either positive definite or negative definite.
\end{itemize}
Given such a matrix $J$, we define an algebraic group $\G$ over $\Q$ by
\[
	\G(A) = \{ g \in \SL_{d+1}(A \otimes_\Q k) : \bar g^t J g = J\},
\]
for any commutative $\Q$-algebra $A$.
We shall regard $\G(\Q)$ as a subgroup of $\G(\R) \times \G(\A_f)$, where $\A_f$ is the ring of finite ad\`eles of $\Q$.
We have an isomorphism $\G(\R) \cong \SU(d,1) \times \SU(d+1)^{e-1}$.
Given any compact open subgroup $K_f\subset \G(\A_f)$, we let $\Gamma(K_f)$ be the group of elements of $\G(\Q)$, which project into $K_f$.
The projection of $\Gamma(K_f)$ in $\SU(d,1)$ is called a
\emph{congruence subgroup of $\SU(d,1)$ of the first kind}.
Any subgroup of $\SU(d,1)$ which is commensurable with $\Gamma(K_f)$ is called an \emph{arithmetic subgroup of the first kind}.

\begin{theorem}[\cite{hill},\cite{stover toledo 2}]
	Let $\Gamma$ be an arithmetic subgroup of $\SU(d,1)$
	of the first kind,
	constructed using a CM field $k$ with $[k:\Q] > 2$.
	Then the groups $\tilde \Gamma$ and $\tilde \Gamma^{(n)}$
	are residually finite.
\end{theorem}

It is not known whether the theorem extends to the case $[k:\Q]=2$,
and this is the question which we investigate in this paper.
The case $[k:\Q]=2$ is geometrically different from the case $e\ge 2$,
since the groups $\Gamma$ are cocompact for $[k:\Q] > 2$ but
have cusps in the case $[k:\Q]=2$.

To describe our result, let $X_\Gamma$ be the locally symmetric space corresponding to $\Gamma$.
We shall write $H^1_!(X_\Gamma,\C)$ for the image in $H^1(X_\Gamma,\C)$
of the cohomology of compact support $H^1_{\compact}(X_\Gamma,\C)$.
We prove the following result:

\begin{theorem}
	\label{thm:main}
	Let $\Gamma$ be an arithmetic subgroup of $\SU(d,1)$
	with cusps (i.e. constructed from an complex quadratic field $k$).
	Assume that there exists an arithmetic subgroup $\Gamma'$
	commensurable with $\Gamma$
	such that $H^1_!(X_{\Gamma'},\C)\ne 0$.
	Then the groups $\tilde \Gamma$ and $\tilde\Gamma^{(n)}$ are all residually finite.
\end{theorem}

We'll briefly discuss the hypothesis that $H^1_!(X_{\Gamma'},\C) \ne 0$.
The cuspidal cohomology is contained in the inner cohomology
(see \cite{harder vol 3}), so if $H^1(X_{\Gamma'},\C)$ contains any non-zero cusp forms then the hypothesis of \autoref{thm:main} holds.
Furthermore, it is known (see for example \cite{shimura}) that there exists a congruence subgroup $\Gamma'$, such that $H^1(X_{\Gamma'},\C)\ne 0$.
However, the hypothesis of \autoref{thm:main} is that
$H^1_!(X_{\Gamma'},\C)\ne 0$ and this is rather stronger.
We shall give an example of a group satisfying this hypothesis at the end of the paper; the author is extremely grateful to Matthew Stover for suggesting this example.

Here are some equivalent formulations of the hypothesis:

\begin{proposition}
	Let $\Gamma$ be a neat arithmetic subgroup of $\SU(d,1)$ of the first kind with cusps.
	Let $\omega \in H^2(X_{\Gamma},\C)$ be the cohomology
	class represented by the invariant K\"ahler form on the symmetric space attached to $\SU(d,1)$.
	Then the following are equivalent:
	\begin{enumerate}
		\item
		$H^1_!(X_\Gamma,\C) \ne 0$;
		\item
		$H^{2d-1}_!(X_\Gamma,\C) \ne 0$;
		\item
		There exists $\phi \in H^1(X_\Gamma,\C)$
		such that $\phi \cup \omega^{d-1} \ne 0$ in $H^{2d-1}(X_\Gamma,\C)$.
	\end{enumerate}
	(Here we are writing $\cup$ for the cup product operation.)
\end{proposition}

The equivalence of 1 and 2 is by duality
(see \eqref{eq:inner pairing}).
The equivalence of 2 and 3 follows immediately from
\autoref{lem:image} below.
In this context, it is reassuring to note that
$\omega^{d-1}$ represents a non-zero cohomology class in $H^{2d-2}(X_\Gamma,\C)$ (see \autoref{lem:invariant cohomology} below).

\bigskip

The paper is organized as follows.
In section 2, we recall a purely group theoretical lemma, which gives a method for showing that certain extension groups are residually finite.
In section 3, we recall some standard facts about the the locally symmetric spaces $X_{\Gamma}$ and their compactifications.
In section 4 we prove \autoref{thm:main}.
In section 5 we give an example of a group $\Gamma$ satisfying $H^1_!(\Gamma,\C) \ne 0$, allowing us to apply \autoref{thm:main} in this case.

\paragraph{Acknowledgements}
The author would like to thank Matthew Stover for suggesting the example in section 5.

\section{A group theoretical lemma}

The method of proof of \autoref{thm:main}
is a modification of the argument in \cite{hill}.
In particular, we shall use the following lemma, which is
proved in both \cite{hill}, \cite{stover toledo} and \cite{stover toledo 2}.
For completeness, we include a short proof.

\begin{lemma}
	\label{lem:criterion}
	Let $G$ be a finitely generated, residually finite group, and suppose that we have a central extension
	\[
		1 \to \Z \to \tilde G \to G \to 1.
	\]
	Let $\sigma_\Z\in H^2(G,\Z)$ be the cohomology class of the extension, and let $\sigma_\C$ be the image of $\sigma_\Z$ in $H^2(G,\C)$.
	Assume that there exist elements $\phi_i,\psi_i\in H^1(G,\C)$
	such that
	\[
		\sigma_\C = \phi_1 \cup \psi_1 + \cdots + \phi_r \cup \psi_r.
	\]
	Then $\tilde G$ is residually finite.
	Furthermore, the quotient group $\tilde G/n\Z$ is residually finite for every positive integer $n$.
\end{lemma}

\begin{proof}
	Let $G^{\ab} = G / [G,G]$.
	Elements of $H^1(G,\C)$ may be regarded as group homomorphisms $G \to \C$.
	Every such homomorphism is the inflation of a homomorphism $G^{\ab} \to \C$.
	Hence $\sigma_\C$ is also the inflation of a cohomology class on $G^{\ab}$, and we shall write $\widetilde{G^{\ab}}$ for the corresponding central extension.
	It follows that we have a commutative diagram with exact rows:
	\begin{center}
	\begin{tikzcd}
		1\ar[r] & \Z\ar[r]\ar[d,hookrightarrow] &
		\tilde G\ar[r]\ar[d] & G\ar[r]\ar[d] & 1\\
		1\ar[r] & \C\ar[r] & \widetilde{G^{\ab}}\ar[r] & G^{\ab}\ar[r] & 1
	\end{tikzcd}
	\end{center}
	We shall write $\Delta$ for the image of $\tilde G$ in $\widetilde{G^{\ab}}$.
	The group $\widetilde{G^{\ab}}$ is nilpotent.
	Hence $\Delta$
	is a finitely generated nilpotent group, and is therefore
	residually finite.
	The group $\tilde G$ injects into $\Delta \times G$.
	Since $\Delta$ and $G$ are both residually finite, it follows that $\tilde G$ is residually finite.
	Similarly, since $\Delta / n\Z$ is residually finite, it follows that $\tilde G / n\Z$ is residually finite.
\end{proof}

Again let $\Gamma$ be an arithmetic subgroup of $\SU(d,1)$ of the first kind with cusps.
\autoref{lem:criterion} will be applied in the case that $G=\Gamma$, $\tilde G$ is its pre-image in $\widetilde{\SU(d,1)}$
and $\tilde G / n\Z= \tilde\Gamma^{(n)}$.
We note that one may construct a different central extension
of $\Gamma$ for which \autoref{lem:criterion} cannot be applied.
For example, take $d=2$ and let $\tau \in H^2(\Gamma,\Z)$ be a
cohomology class whose restriction to the Borel--Serre boundary of $X_\Gamma$ is non-torsion (for example an Eisenstein cohomology class, see \cite{harder}).
Such a class $\tau$ cannot be
expressed as a sum of cup products of elements of $H^1$,
because all such cup products restrict to torsion on the Borel--Serre boundary.

\section{Background material}

We shall now recall the construction of arithmetic subgroups of $\SU(d,1)$ with cusps.
Let $k$ be a complex quadratic extension of $\Q$; we shall
identify $k$ with a subfield of $\C$, and we shall write $z \mapsto \bar z$ for complex conjugation on $\C$ or on $k$.
Let $J_0$ be a $(d-1)\times (d-1)$ positive definite Hermitian matrix with entries in $k$
and let
\[
	J = \begin{pmatrix}
		0 & 0 & 1 \\
		0 & J_0 & 0 \\
		1 & 0 & 0
	\end{pmatrix}.
\]
The matrix $J$ defines a Hermitian form on $\C^{n+1}$ of signature $(d,1)$ by
\[
	( v, w ) = \bar v^t J w,
\]
where $\bar v^t$ denotes the conjugate transpose of a column matrix $v$.

We define an algebraic group $\G$ over $\Q$ to be the group of
isometries in $\SL_{d+1}$ of the Hermitian form.
More precisely, for a $\Q$-algebra $A$, we define
\[
	\G( A) = \{ g \in \SL_3(A \otimes_\Q k) :
	\bar g^t J g = J \}.
\]
Since the matrix $J$ has signature $(d,1)$, the group $\G(\R)$
may be identified with $\SU(d,1)$.

Let $\A_f$ be the ring of finite ad\`eles of $\Q$.
The group $\G(\A_f)$ is totally disconnected, and contains the projection of $\G(\Q)$
as a dense subgroup (by Kneser's Strong Approximation Theorem).
For a compact open subgroup $K_f \subset \G(\A_f)$, the intersection $\Gamma(K_f) = \G(\Q) \cap K_f$ is called a \emph{congruence subgroup} of $\G(\Q)$.
Any subgroup of $\G(\Q)$ which is commensurable with a congruence subgroup is called an \emph{arithmetic subgroup}.
It is known that there exist arithmetic subgroups of $\G(\Q)$, which are not congruence subgroups.

The group $\G(\A_f)$ may be identified with the projective limit
of the sets $\G(\Q)/\Gamma(K_f)$, where $\Gamma(K_f)$ ranges over the congruence subgroups of $\G(\Q)$.
We also define the \emph{arithmetic completion} $\widehat{\G(\Q)}$
to be the projective limit of the sets $\G(\Q)/\Gamma$, where $\Gamma$ ranges over all the arithmetic subgroups of $\G(\Q)$.
Since the filtration by arithmetic subgroups is invariant under conjugation, the arithmetic completion is a group, and we have a
natural surjective homomorphism $\widehat{\G(\Q)} \to \G(\A_f)$.
The kernel of the homomorphism is an infinite profinite group, and is called the \emph{congruence kernel} $C_\G$ of $\G$.

An arithmetic group $\Gamma$ is said to be \emph{neat} if for every $g \in \Gamma$, the eigenvalues of $g$ generate a torsion-free subgroup of $\C^\times$.
For every congruence subgroup $\Gamma$, there is a neat congruence subgroup $\Gamma'$ of finite index in $\Gamma$.
Every subgroup of a neat group is neat, and every neat group is torsion-free.

\subsection{Quotient spaces and compactifications}

Let
\[
	\cH = \left\{ [v] \in \Proj^d(\C) :
	( v,v ) < 0 \right\},
\]
where we are writing $[v]$ for the point in projective space
represented by a non-zero vector $v$.
The complex manifold $\cH$ has an obvious action of $\SUd$,
 and is a model of the symmetric space attached to $\SU(d,1)$.
For each arithmetic subgroup $\Gamma$ in $\G(\Q)$, we shall write $X=X_\Gamma$ for the quotient space $\Gamma\backslash \cH$.
If $\Gamma$ is neat then $X_\Gamma$ is a smooth, non-compact complex manifold.

By a \emph{cusp}, we shall mean a point $[v]$ of $\Proj^d(k)$, such that
$( v,v) =0$.
For each such $[v]$, there is a parabolic subgroup
$P_v$ of $\SU(d,1)$, defined by
\[
	P_v
	= \{ g \in \SU(d,1) :  [g \cdot v] = [v]\}.
\]
If $\Gamma$ is an arithmetic subgroup of $\G(\Q)$, then $\Gamma$ permutes the cusps with finitely many orbits.

Assume that $[v]$ is a cusp, with
corresponding parabolic subgroup $P_v$.
We may choose a Langlands decomposition
\[
	P_v = M_v A_v N_n,
\]
where $A_v$ is the connected component of a split torus in $P_v$ which is isomorphic to $\R^{>0}$;
the group $M_v$ is isomorphic to $\U(d-1)$,
and $N_v$ is the unipotent radical of $P_v$.
There is a homomorphism $\phi_v : P_v \to \R^{>0}$ defined by $\phi_v(p) = |\lambda|$, where $\lambda\in\C$ satisfies $p v =\lambda \cdot v$.
The subgroups $M_v$ and $N_n$ are in the kernel of $\phi_v$, and the restriction of $\phi_v$ to $A_v$ is an isomorphism.

By the Iwasawa decomposition,
the group $A_v \ltimes N_v$ acts simply transitively
on the symmetric space $\cH$, so by choosing a base point, we may identify $\cH$ with this group.

Assume from now on that $\Gamma$ is a neat arithmetic subgroup of $\G(\Q)$.
For such groups $\Gamma$,
the intersection $\Gamma_v=\Gamma\cap P_v$
is contained in $N_v$, and is a cocompact
subgroup of $N_v$.
The subgroup $\Gamma_v$
acts on $A_v \ltimes N_v$ by translation on $N_v$, preserving the $A_v$-coordinate.
It therefore acts also on the following subset
\[
	A_v^{< \epsilon} N_v
	= \{ a n : a\in A_v, n\in N_v, |\phi_v(a)| < \epsilon \},
\]
where $\epsilon$ is a positive real number.
We may choose $\epsilon$ sufficiently small so that the quotient space $U_v = \Gamma_v \backslash (A_v^{< \epsilon} N_v)$
injects into $X_\Gamma$.
We shall call such a subset $U_v$ of $X_\Gamma$ a \emph{neighbourhood of the cusp $[v]$}.
By reduction theory, there are finitely many
non-intersecting cusp neighbourhoods (one for each $\Gamma$-orbit of cusps), such that the complement of the cusp neighbourhoods is a
compact subset of $X_\Gamma$.

As an example, choose the cusp $v = \begin{pmatrix}
1\\0 \\ \vdots \\ 0 
\end{pmatrix}$.
In this case $N_v$ is a Heisenberg group, consisting
of all matrices of the form
\[
	n(z,x)
	=
	\begin{pmatrix}
		1 & -\bar z^t J_0 & -\frac{|| z||^2}{2} + i x \\
		0 & I_{n-1} & z\\
		0 & 0 & 1
	\end{pmatrix},
	\qquad
	z \in \C^{d-1},
	x\in \R,
\]
where $|| z||^2 = \bar z^t J_0 z$.
We have a short exact sequence of Lie groups:
\[
	\begin{matrix}
	1 &\to& \R &\to &N_v& \to &\C^{d-1}& \to& 1\\
	&& x &\mapsto & n(0,x)\\
	&&&& n(z,x) & \mapsto & z
	\end{matrix}
\]
The image of $\Gamma_v$ in $\C^{d-1}$ is a full lattice $L_v$,
and the quotient $\C^{d-1}/L_v$ is an abelian variety
(indeed $L_v$ is commensurable with $\cO_k^{d-1}$).
The kernel of the map
$\Gamma_v \to L_v$ is isomorphic to $\Z$.
Hence the topological space $\Gamma_v \backslash N_v$ is an $\R/\Z$ bundle over the abelian variety $\C^{d-1}/L_v$.
The cusp neighbourhood $U_v$ is the product of this space with the open interval $(0,\epsilon)$.

If we choose any other cusp $v$, then
the groups $P_v$ and $N_v$ are conjuagates in $\SU(d,1)$ of the subgroups described above.
Hence the cusp neighbourhood has a similar description as a product of an open interval with
a circle bundle over an abelian variety.

We shall consider two compactifications of $X$.
The first is the Borel--Serre compactification
in which we embed each cusp neighbourhood $U_v$
into a larger topological space $U_v^{\BS}$ as follows:
\begin{center}
	\begin{tikzcd}
		U_v \ar[r,equal] \ar[d,phantom,sloped,"\subset"] &  \Gamma_v \backslash  N_v \times (0,\epsilon) \\
		U_v^{\BS} \ar[r,equal]&
		\Gamma_v \backslash  N_v \times [0,\epsilon) .
	\end{tikzcd}
\end{center}
The embedding $U_v \to U_v^{\BS}$
is evidently a homotopy equivalence.
Therefore the resulting compactification
$X^{\BS}$ has the same cohomology groups as $X$.

We shall write $\partial X^{\BS}$ for the complement of $X$ in its Borel--Serre compactification.
The boundary of the Borel--Serre compactification is a disjoint union of
manifolds homeomorphic to
\[
	\partial X_{v}^{\BS}
	= \Gamma_v \backslash N_v.
\]
Each of the boundary components $\partial X_{v}^{\BS}$
is a circle bundle over an abelian variety.

The second compactification which we shall consider is the smooth compactification $\tilde X$ constructed in \cite{AshMumfordRapoportTai}.
This may be obtained from the Borel--Serre compactification by quotienting each boundary component $\Gamma_v \backslash N_v$ by
the centre of $N_v$, i.e. by the subgroup
of matrices of the form $n(0,x)$.
The resulting boundary component is
the abelian variety
\[
	\partial \tilde X_v = \C^{d-1}/L_v.
\]
The compactification $\tilde X$ is a smooth complex manifold
but is not homotopic to to $X$.
There is an obvious projection map
$X^{\BS} \to \tilde X$.

We shall write $\partial \tilde X$ for the complement of $X$ in its smooth compactification.
This boundary set $\partial \tilde X$ is the disjoint union of
the abelian varieties $\partial\tilde X_v$.

\subsection{Cohomology groups}

\label{sec:cohomology}

In this paper we shall use various cohomology groups.
For convenience, we list the notation and some standard properties for each of these.
In almost all cases we shall consider cohomology with coefficients in $\C$. In such cases, we shall not always write in the coefficients.

\begin{itemize}
	\item
	The continuous cohomology groups of the group $\SUd$
	will be written $H^\bullet_{\cts}(\SUd,-)$.
	We may identify $H^r_{\cts}(\SUd,\C)$ with the vector space of differential $r$-forms on $\cH$, which are invariant under the action of $\SU(d,1)$
	(see \cite{BorelWallach}).

	For example, there is an invariant K\"ahler form
	$\omega$ on $\cH$.
	The form $\omega$ generates $H^2_{\cts}(\SU(d,1),\C)$.
	More generally we have
	\[
		H^r_{\cts}(\SU(1,d),\C)
		=\begin{cases}
			\C \cdot \omega^{r/2} & \text{if $r=0,2,4,\ldots,2d$}\\
			0 & \text{otherwise.}
		\end{cases}
	\]
	\item
	The measurable cohomology groups
	(see \cite{moore}) of a connected Lie group $G$ will be written $H^\bullet_{\meas}(G,-)$.
	For a connected Lie group $G$ with fundamental group $\pi_1(G)$, there is a canonical isomorphism
	$H^2_{\meas}(G,\Z) \cong \Hom(\pi_1(G),\Z)$.
	
	In particular, the group $\SU(d,1)$ has fundamental group $\Z$, so we have $H^2_{\meas}(\SU(d,1),\Z) \cong \Z$.
	We shall choose a generator
	$\sigma_\Z$ for this group, i.e.
	\[
		H^2_{\meas}(\SU(d,1),\Z) = \Z\cdot  \sigma_\Z.
	\]
	The group extension of $\SU(d,1)$ corresponding to the cocycle
	$\sigma_\Z$ is the universal cover of $\SU(d,1)$.
	
	By \cite{wigner} there is an isomorphism
	\begin{equation}
		\label{eqn:wigner}
		H^\bullet_{\cts}(G,\C)
		\cong
		H^\bullet_{\meas}(G,\Z) \otimes \C.
	\end{equation}
	In particular, the image of $\sigma_\Z$ in $H^2_{\cts}(\SU(d,1),\C)$ is a non-zero multiple of $\omega$.

	\item
	For an arithmetic subgroup $\Gamma$, the Eilenberg--MacLane cohomology groups will be written $H^\bullet(\Gamma,-)$.
	There are restriction maps $H^r_{\cts}(\SU(d,1),\C) \to H^r(\Gamma,\C)$.
	Some of these maps are injective and others are zero.
	(In fact, we'll see in \autoref{lem:invariant cohomology} that the map is injective if $r < 2d$
	and zero if $r=2d$).

	\item
	If $\Gamma$ is a neat arithmetic subgroup
	of $\SUd$, then the quotient space
	$X= \Gamma \backslash \cH$ is a complex $d$-dimensional manifold.
	We shall write $H^\bullet(X)$ for the singular or de Rham cohomology groups of this manifold with complex coefficients.
	
	Apart from the manifold $X$, we shall also consider two compactifications $X^{\BS}$ and $\tilde X$.
	
	Recall that there are canonical isomorphisms
	\[
		H^\bullet(X^{\BS}) \cong H^\bullet(X) \cong H^\bullet(\Gamma,\C).
	\]
	The composition
	\[
		H^\bullet_{\cts}(\SU(d,1),\C)
		\stackrel{\Rest}\to H^\bullet(\Gamma,\C)
		\cong H^\bullet(X),
	\]
	takes an invariant differential form on $\cH$ to its de Rham cohomology class on $X$.

	\item
	We shall write $c_1(X)$ for the first Chern class of the
	canonical sheaf on $X$, regarded as an element of $H^2(X)$.
	It is known that $c_1(X)$ is a multiple of the cohomology class
	$\omega$ by a positive real number.

	\item
	We shall write $H^\bullet_{\compact}(X)$ for the
	compactly supported cohomology of $X$ with complex coefficients.
	The space $H^{2d}_{\compact}(X)$ is one-dimensional, and the cup-product map
	\[
		\cup : H^r(X) \otimes H^{2d-r}_{\compact}(X) \to H^{2d}_{\compact}(X) \cong \C
	\]
	is a perfect pairing, allowing us to identify $H^r(X)$ with the dual space of $H^{2d-r}_{\compact}(X)$.
	
	\item
	We shall write $\partial X^{\BS}$ and $\partial \tilde X$ for the boundaries of the compactifications.
	The relative cohomology groups
	$H^\bullet(X^{\BS},\partial X^{\BS})$ and
	$H^\bullet(\tilde X,\partial \tilde X)$
	are both canonically isomorphic to
	the compactly supported cohomology.
	We therefore have a commutative diagram whose
	rows are long exact sequences:
	\begin{center}
		\begin{tikzcd}
			\cdots\ar[r] &
			H^{\bullet}_{\compact}(X) \ar[r]\ar[d,equal]  &
			H^\bullet(\tilde X) \ar[r] \ar[d]&
			H^\bullet(\partial \tilde X) \ar[r]\ar[d] &
			H^{\bullet+1}_{\compact}(X) \ar[r]\ar[d,equal]& \cdots \\
			\cdots\ar[r]
			& H^{\bullet}_{\compact}(X) \ar[r]
			& H^\bullet(X)  \ar[r]
			& H^\bullet(\partial X^{\BS}) \ar[r]
			& H^{\bullet+1}_{\compact}(X) \ar[r]
			&\cdots
		\end{tikzcd}
	\end{center}
	We shall write $H^\bullet_!(X)$ for the kernel of the restriction map
	$H^\bullet(X) \to H^\bullet(\partial X^{\BS})$,
	or equivalently the image of the map
	$H^\bullet_{\compact}(X) \to H^\bullet(X)$.
	The vector spaces $H^\bullet_!(X)$ are known as the \emph{inner cohomology groups}.
	\item
	For $r=0,\ldots, 2d$ there is a perfect pairing (see \cite{harder vol 3}):
	\begin{equation}
		\label{eq:inner pairing}
		\langle -,-\rangle : H^r_!(X) \otimes H^{2d-r}_!(X)
		\to H^{2d}_{\compact}(X) \cong \C,
	\end{equation}
	defined as follows. Given $a\in H^r_!(X)$
	and $b\in H^{2d-r}_!(X)$, we may choose a pre-image $a_{\compact} \in H^r_{\compact}(X)$ of $a$.
	The pairing $\langle a,b \rangle$ is defined to be the cup product $a_{\compact} \cup b$.
	This cup product does not depend on the choice of $a_{\compact}$.
	\item
	We shall use the notation
	\[
		H^\bullet_{\stable}
		=
		\lim_{\genfrac{}{}{0pt}{1}{\longrightarrow}{\Gamma'}}
		H^\bullet(X_{\Gamma'}),
		\qquad
		H^\bullet_{!,\stable}
		=
		\lim_{\genfrac{}{}{0pt}{1}{\longrightarrow}{\Gamma'}}
		H^\bullet_!(X_{\Gamma'}),
		\qquad
		H^\bullet_{\compact,\stable}
		=
		\lim_{\genfrac{}{}{0pt}{1}{\longrightarrow}{\Gamma'}}
		H^\bullet_{\compact}(X_{\Gamma'}),
	\]
	where the limits are taken over all arithmetic subgroups $\Gamma'$ of $\Gamma$.
	These direct limits may be regarded as unions, since all of the connecting homomorphisms are injective.

	There is an obvious action of $\G(\Q)$ on the vector spaces
	$H^\bullet_{\stable}$ and $H^\bullet_{!,\stable}$, and this action extends to a smooth action of the totally disconnected group  $\widehat{\G(\Q)}$.

	\item
	Since cup products are compatible with restriction maps, the pairing (\ref{eq:inner pairing}) extends to
	a perfect pairing
	\[
		H^r_{!,\stable} \otimes H^{2d-r}_{!,\stable} \to H^{2d}_{\compact,\stable} \cong \C.
	\]
	This pairing is $\widehat{\G(\Q)}$-invariant, in the sense that
	\[
		\langle g a , g b \rangle
		=
		\langle a , b \rangle
		\quad
		\text{for all $g \in \widehat{\G(\Q)}$.}
	\]
	\item
	There is also an invariant positive definite inner product
	on $H^r_{!,\stable}$ defined by
	\[
		\langle \langle a,b\rangle \rangle = \langle a, *b \rangle,
	\]
	(see page 73 of \cite{harder vol 3}).
	In particular, the representation $H^r_{!,\stable}$ of $\widehat{\G(\Q)}$ is semi-simple.

	\item
	We shall use the notation
	\[
		H^\bullet_{!,\Cong}
		=
		\lim_{\genfrac{}{}{0pt}{1}{\longrightarrow}{K_f}}
		H^\bullet_!(X_{\Gamma(K_f)}),
	\]
	where the limit is taken over all congruence subgroups $\Gamma(K_f)$.	
	There is a smooth action of $\G(\A_f)$ on the vector space
	$H^\bullet_{!,\Cong}$, and we may identify $H^\bullet_{!,\Cong}$
	with the subspace of invariants $(H^\bullet_{!,\stable})^{C_\G}$,
	where $C_\G$ is the congruence kernel of $\G$.

	\item
	The vector space $H^\bullet_{!,\Cong}$ is
	a semi-simple representation of $\G(\A_f)$.
	More precisely there is a countable direct sum decomposition
	(see \cite{harder vol 3}):
	\begin{equation}
		\label{eqn:matsushima}
		H^r_{!,\Cong}
		=
		\bigoplus_{\pi \in \Pi_!^r}
		H^r_{\cts}(\SUd,\pi_\infty) \otimes \pi_f,
	\end{equation}
	where $\Pi_!^r$ is a certain set of automorphic representations of $\G(\A)$.
	Each of the automorphic representations $\pi$ decomposes as
	$\pi_\infty \otimes \pi_f$, where $\pi_\infty$ is a simple representation of $\SUd$
	and $\pi_f$ is a smooth, simple representation of $\G(\A_f)$.
	\item
	We shall be particularly interested in the subspace
	$(H^\bullet_{!,\stable})^{\widehat{\G(\Q)}}$ of $\widehat{\G(\Q)}$-invariant cohomology classes.
	Since $(H^\bullet_{!,\stable})^{C_\G}= H^\bullet_{!,\Cong}$, it follows that
	\begin{equation}
		\label{eqn:stable invariants}
		(H^\bullet_{!,\stable})^{\widehat{\G(\Q)}}
		=
		(H^\bullet_{!,\Cong})^{\G(\A_f)}.	
	\end{equation}
	The right hand side of \eqref{eqn:stable invariants} may be evaluated using \eqref{eqn:matsushima}.
	The trivial representation $\C$ occurs with multiplicity at most $1$ in each set $\Pi^r_!$.
	For all non-trivial representations $\pi$ in $\Pi^r_!$, the vector space $\pi_f$ is infinite dimensional.
	This implies
	\begin{equation}
		\label{eq:invariant inner}
		\left(H^r_{!,\stable}\right)^{\widehat{\G(\Q)}}
		\cong
		\begin{cases}
			H^r_{\cts}(\SU(d,1),\C) &
			\text{if the trivial representation is in $\Pi_!^r$,}\\
			0 & \text{otherwise.}
		\end{cases}
	\end{equation}
	By the discussion above, if $\left(H^r_{!,\stable}\right)^{\widehat{\G(\Q)}}$ is non-zero, then $r$ is even and this space is spanned by the
	cohomology class of $\omega^{r/2}$ on $X_\Gamma$ (or more accurately, by the image of this cohomology class in the direct limit $H^r_{\stable}$).

\end{itemize}

\begin{lemma}
	\label{lem:omega non-zero}
	Let $\Gamma$ be an arithmetic subgroup of $\G(\Q)$.
	The invariant  K\"ahler form $\omega$
	represents a non-zero cohomology class on $X_\Gamma$.
	Consequently, the restriction
	of $\sigma_\Z$ to $\Gamma$ is a non-zero element of $H^2(\Gamma,\Z)$.
\end{lemma}

(Here we are using the assumption that $d \ge 2$; the statement would be false for $\SU(1,1)$.)

\begin{proof}
	Since $\omega$ is a multiple of the image of $\sigma_\Z$,
	it is sufficient to prove the statement for $\omega$.
	It is also sufficient to prove the result with $\Gamma$ sufficiently small.
	We may therefore assume, without loss of generality,
	that $\Gamma$ is neat.
	Hence $X$ is a smooth complex manifold.

	Suppose for a moment that there exists a compact Riemann surface $Y \subset X$.
	The K\"ahler form $\omega$ restricts to a K\"ahler form on $Y$, and by positivity of K\"ahler forms (see for example \cite{GriffithsHarris}) we have
	\[
		\int_{[Y]} \omega > 0.
	\]
	Therefore $\omega$ represents a non-zero cohomology class on $Y$, and hence also on $X$.

	It is therefore sufficient to find a compact Riemann surface $Y$ contained in $X$.
	The folloing construction of such a $Y$ is taken from \cite{shimura}.		
	We may choose a 2-dimensional subspace $S \subset k^{d+1}$, such that
	the Hermitian form has signature $(1,1)$ on $S$ and is anisotropic.
	Let $G$ be the group of isometries of $S \otimes_{\Q} \R$.
	Our choice of $S$ implies that $G$ is isomorphic to $\SU(1,1)$
	and $\Gamma \cap G$ is cocompact in $G$.
	Let $Y$ be the locally symmetric space corresponding to the
	subgroup $\Gamma \cap G$ of $G$.
	The inclusion of $G$ in $\SU(d,1)$ gives us an inclusion of
	$Y$ in $X$ as a compact Riemann surface.
\end{proof}

\begin{lemma}
	\label{lem:Gamma v restriction}
	Let $\Gamma$ be a neat arithmetic subgroup of $\G(\Q)$.
	For any cusp $[v]$, the restrictions of
	$\sigma_\Z$ and $\omega$ to $\Gamma_v$ are coboundaries.
\end{lemma}

\begin{proof}
	As $\omega$ is a multiple of the image of $\sigma_\Z$ in
	$H^2_{\cts}(\SU(d,1),\C)$, it is sufficient to prove the result for $\sigma_\Z$.
	Since $\Gamma$ is neat, we have $\Gamma_v \subset N_v$.
	Therefore the restriction map $H^2_{\meas}(\SU(d,1),\Z) \to H^2(\Gamma_v,\Z)$ factors through the group
	$H^2_{\meas}(N_v,\Z)$.
	Since the Lie group $N_v$ is simply connected, we have $H^2_{\meas}(N_v,\Z)=0$.
	Therefore the restriction of $\sigma_\Z$ to
	$\Gamma_v$ is zero.
\end{proof}

\begin{lemma}
	\label{lem:omega is inner}
	Let $\Gamma$ be a neat arithmetic subgroup
	of $\G(\Q)$ and let $X=X_{\Gamma}$.
	The image of $\omega$ in $H^2(X)$ is in
	the subspace $H^2_!(X)$ of inner cohomology classes.
\end{lemma}

\begin{proof}
	It is sufficient to show that the class $\omega$ vanishes on each Borel-Serre boundary component $\partial X^{\BS}_v$.
	This follows from \autoref{lem:Gamma v restriction}, in view of the following commutative diagram.
	\begin{center}
		\begin{tikzcd}
			H^2(X^{\BS})\ar[r]\ar[d,phantom,sloped,"\cong"] &H^2(\partial X^{\BS}_v) \ar[d,phantom,sloped,"\cong"] \\
			H^2(\Gamma,\C)\ar[r]  &H^2(\Gamma_v,\C)
		\end{tikzcd}
	\end{center}
\end{proof}

\begin{lemma}
	\label{lem:sigma stable}
	The vector space $(H^{2}_{!,\stable})^{\widehat{\G(\Q)}}$ is
	one-dimensional, and is spanned by
	the cohomology class of the invariant K\"ahler form
	$\omega$ on $X_\Gamma$.
\end{lemma}

\begin{proof}
	By \eqref{eq:invariant inner}, the space $(H^2_{!,\stable})^{\widehat{\G(\Q)}}$ is at most one-dimensional, so it is sufficient to show that $\omega$ is a non-zero element in this space.
	\autoref{lem:omega is inner} shows that $\omega$ is in the subspace of
	inner cohomology classes, and
	\autoref{lem:omega non-zero} implies that
	$\omega$ represents a non-zero cohomology class on $X_\Gamma$ for every arithmetic subgroup $\Gamma$.
\end{proof}

\begin{lemma}
	\label{lem:invariant cohomology}
	Let $\Gamma$ be an arithmetic subgroup of $\G(\Q)$.
	For $r=1,\ldots,d-1$, the $2r$-form $\omega^r$
	represents a non-zero cohomology class on $X_\Gamma$,
	which spans $(H^{2r}_{!,\stable})^{\widehat{\G(\Q)}}$.
\end{lemma}

\begin{proof}
	By \autoref{lem:sigma stable}, the representation $H^2_{!,\stable}$ of $\widehat{\G(\Q)}$
	has a trivial 1-dimensional subrepresentation spanned by $\omega$.
	Hence by duality, $H^{2d-2}_{!,\stable}$ has a trivial
	1-dimensional quotient representation.
	By semi-simplicity, it follows that $(H^{2d-2}_{!,\stable})^{\widehat{\G(\Q)}}$ is non-zero.
	By the discussion above, $(H^{2d-2}_{!,\stable})^{\widehat{\G(\Q)}}$ is one-dimensional, and is spanned by
	$\omega^{d-1}$.
	In particular, $\omega^{d-1}$ represents a non-zero cohomology class on $X_\Gamma$.
	From this, it follows that $\omega^{r}$ represents a non-zero cohomology class on $X_\Gamma$ for $1\le r \le d-1$.
	The image of $\omega^r$ in $H^{2r}_{!,\stable}$ spans a trivial one-dimensional subrepresentation.
	Therefore $(H^{2r}_{!,\stable})^{\widehat{\G(\Q)}}\ne 0$.
	By the discussion above, $(H^{2r}_{!,\stable})^{\widehat{\G(\Q)}}$ is spanned by $\omega^r$.
\end{proof}

\begin{proposition}
	\label{prop:sigma extend}
	Let $\Gamma$ be a neat arithmetic subgroup of $\G(\Q)$
	and let $X=X_\Gamma$.
	Then there exists an element $\tilde \omega \in H^{2}(\tilde X)$, such that
	\begin{itemize}
		\item
		the restriction of $\tilde \omega$ to $X$ is
		the invariant K\"ahler form $\omega$.
		\item
		$\tilde\omega$ is in the ample cone in $H^{1,1}(\tilde X)$.
	\end{itemize}
\end{proposition}

\begin{proof}
	The lemma does not depend on our choice of normalization of
	$\omega$; we shall assume for simplicity that $\omega = c_1(X)$ in $H^2(X)$.
	
	Let $\tilde \omega = c_1(\tilde X) + \epsilon \cdot [\partial \tilde X]$, where $c_1(\tilde X)$ is the first Chern class of the canonical sheaf on $\tilde X$, and
	$\epsilon$ is a positive real number.
	Here we are writing $[\partial\tilde X]$ for the Poincar\'e dual of the $2d-2$-cycle $\partial \tilde X$, or equivalently the first Chern class of the line bundle corresponding the the divisor $\partial \tilde X$.
	It is shown in \cite{di cerbo di cerbo}, that if $\epsilon$ is
	sufficiently small then $\tilde \omega$ is in the ample cone.
	By naturality of Chern classes, it follows that the restriction of $c_1(\tilde X)$ to $X$ is $c_1(X)$, which we are assuming is equal to $\omega$.
	The restriction of the divisor $\partial \tilde X$ to $X$ is $0$; hence the restriction of $[\partial \tilde X]$ to $X$ is $0$.
	It follows that the restriction of $\tilde\omega$ to $X$ is
	$\omega$.
\end{proof}

\section{Proof of Theorem 1}

Fix a neat arithmetic subgroup $\Gamma \subset \G(\Q)$ and let $X$ be the quotient space $\Gamma\backslash \cH$.
We shall let $\tilde X$ denote the smooth compactification of $X$.
We shall write $\partial \tilde X$ for the union of the boundary components of $\tilde X$. Each boundary component is an abelian variety.
We choose a neighbourhood $\tilde U$ of $\partial\tilde X$, so that $\partial \tilde X$ is a deformation retract of $\tilde U$.
We shall also write $U$ for the intersection $\tilde U\cap X$.
Note that $U$ is homotopic to the Borel--Serre boundary of $X$.
The Mayer--Vietoris sequence for the cover $\tilde X = X \cup \tilde U$ takes the form:
\begin{equation}
	\label{eqn:Mayer--Vietoris}
	\to H^n(\tilde X) \to H^n(X) \oplus H^n(\partial \tilde X)
	\to H^n(\partial X^{\BS})
	\to H^{n+1}(\tilde X) \to .
\end{equation}

Let $[v]\in \Proj^d(k)$ be a cusp, and let $N_v$ be the unipotent radical in the parabolic subgroup fixing $[v]$.
Recall that we have a central group extension:
\[
	1 \to \R \to N_v \to \C^{d-1} \to 1.
\]
We shall write $\Gamma_v$ for the intersection of $\Gamma$ with $N_v$.
We also let $L_v$ be the image of $\Gamma_v$ in $\C^{d-1}$ and
$Z_v$ be the kernel of the map $\Gamma_v \to L_v$.

The next three lemmas are well known (for example, see \cite{harder}).
We include proofs for the sake of completeness.

\begin{lemma}
	\label{lem:restriction zero}
	The restriction map $H^1(\Gamma_v,\C) \to H^1(Z_v,\C)$ is zero.
\end{lemma}

\begin{proof}
	We shall regard elements of $H^1(\Gamma,\C)$ as group homomorphisms
	$\phi : \Gamma_v \to \C$.
	We must prove that $\phi(g)=0$ for all elements $g \in Z_v$.
	For any such $g$, there is a positive integer $n$ such that $g^n \in [\Gamma_v,\Gamma_v]$.
	Therefore $\phi(g)= \frac{1}{n}\phi(g^n) = 0$.  
\end{proof}

\begin{lemma}
	\label{lem:H1 boundary comparison}
	The pullback map $H^1(\partial\tilde X) \to H^1(\partial X^{\BS})$ is an isomorphism.
\end{lemma}

\begin{proof}
It is sufficient to show that for each cusp $v$,
the pullback $H^1(\partial\tilde X_v) \to H^1(\partial X_v^{\BS})$ is an isomorphism.
Recall that $\partial\tilde X_v$ is an abelian variety $\C^{d-1}/L_v$, and $\partial X_v^{\BS}$
is a circle bundle over this abelian variety,
homeomorphic to $\Gamma_v\backslash N_v$.
We shall write $Z_v$ for the kernel of the homomorphism $\Gamma_v \to L_v$.
We therefore have a commutative diagram
\begin{center}
	\begin{tikzcd}
		& H^1(\partial \tilde X_v) \ar[r]\ar[d,phantom,sloped,"\cong"]& H^1(\partial X_v^{\BS}) \ar[d,phantom,sloped,"\cong"]\\
		0 \ar[r] & H^1(L_v,\C) \ar[r] & H^1(\Gamma_v,\C) \ar[r] &H^1(Z_v,\C)\;,
	\end{tikzcd}
\end{center}
where the bottom row is the inflation--restriction sequence in group cohomology.
The result now follows from \autoref{lem:restriction zero}.
\end{proof}

\begin{lemma}
	\label{lem:techincal}
	The restriction map gives an isomorphism
	$H^1(\tilde X) \cong H^1(X)$.
\end{lemma}

\begin{proof}
	Consider the following section of the Mayer-Vietoris sequence \eqref{eqn:Mayer--Vietoris}:
	\[
		H^0(X) \oplus H^0(\partial \tilde X)
		\to H^0(\partial X^{\BS}) \to H^1(\tilde X) \to H^1(X) \oplus H^1(\partial \tilde X)
		\to H^1(\partial X^{\BS}).
	\]
	The map $H^0(\partial \tilde X)	\to H^0(\partial X^{\BS})$ is
	clearly an isomorphism.
	By \autoref{lem:H1 boundary comparison}, the pull-back map $H^1(\partial \tilde X)\to H^1(\partial X^{\BS})$ is an isomorphism.
	Hence $H^1(\tilde X) \cong H^1(X)$.
\end{proof}

\begin{lemma}
	\label{lem:image}
	The map $H^1(X) \to H^{2d-1}(X)$ given by
	cup product with $\omega^{d-1}$ has image $H^{2d-1}_!(X)$.
\end{lemma}

\begin{proof}
By \autoref{prop:sigma extend}, there exists
an ample class $\tilde \omega \in H^2(\tilde X)$,
whose restriction to $X$ is the class $\omega$.
We have a commutative diagram:
\begin{center}
	\begin{tikzcd}
		H^1(\tilde X)  \ar[r,"\cong"] \ar[d,"\tilde\omega^{d-1}"]
		& H^1(X) \ar[d,"\omega^{d-1}"]\\
		H^{2d-1}(\tilde X) \ar[r]
		& H^{2d-1}(X) \ar[r]
		& H^{2d-1}(\partial X^{\BS})\; ,
	\end{tikzcd} 
\end{center}
where the vertical maps are given by cup product with $\tilde\omega^{d-1}$ and $\omega^{d-1}$ respectively.
The bottom row is exact, as it is part of the
Mayer--Vietoris sequence (\ref{eqn:Mayer--Vietoris}). Note that $H^{2d-1}(\partial \tilde X)=0$ because $\partial \tilde X$ has dimension $2d-2$.

The diagram is commutative because the restriction of $\tilde{\omega}$ to $X$ is $\omega$.
Since $\tilde{\omega}$ is ample, the Hard Lefschetz Theorem
implies that the left hand vertical map is an isomorphism.
Therefore the image of the right hand vertical map is
$H^{2d-1}_!(X)$.
\end{proof}

\begin{lemma}
	\label{lem:existence}
	Assume $H^1_!(X)\ne 0$.
	Then there exist $\phi \in H^1_!(X)$ and $\psi\in H^1(X)$
	such that $\langle \phi \cup\psi , \omega^{d-1}\rangle \ne 0$.
\end{lemma}

\begin{proof}
	Choose a non-zero $\phi\in H^1_!(X)$.
	Then there exists an element $\phi^* \in H^{2d-1}_!(X)$,
	such that $\langle \phi , \phi^*\rangle \ne 0$.
	By \autoref{lem:image}, we have $\phi^* = \psi\cup\omega^{d-1}$ for some $\psi \in H^1(X)$.
	
	By \autoref{lem:omega is inner} we have $\omega \in H^2_!(X)$,
	so we may choose a pre-image $\omega_{\compact}$ of $\omega$ in $H^2_{\compact}(X)$.
	By definition of the pairing  $\langle -,-\rangle$ in (\ref{eq:inner pairing}), we have
	\begin{align*}
		\langle \phi\cup\psi,\omega^{d-1} \rangle
		&=(\phi \cup \psi) \cup \omega_{\compact}^{d-1} \\
		&=\phi \cup (\psi \cup \omega_{\compact}^{d-1}) \\
		&=\langle \phi, \psi \cup \omega^{d-1}\rangle \\
		&=\langle \phi, \phi^*\rangle  \\
		&\ne 0.
	\end{align*}
	In the third equality above, we have used the fact that
	$\psi\cup \omega_{\compact}^{d-1} \in H^{2d-1}_{\compact}(X)$
	is a pre-image of the element $\psi\cup \omega^{d-1} \in H^{2d-1}_!(X)$.
\end{proof}

\begin{thm}
	Let $\Gamma$ be an arithmetic subgroup of $\SU(d,1)$
	with cusps (i.e. constructed from an complex quadratic field $k$).
	Assume that there exists an arithmetic subgroup $\Gamma'$
	commensurable with $\Gamma$
	such that $H^1_!(X_{\Gamma'})\ne 0$.
	Then the groups $\tilde \Gamma$ and $\tilde\Gamma^{(n)}$ are all residually finite.
\end{thm}

\begin{proof}
	Consider the subspace $V$ of $H^2_{!,\stable}$
	spanned by cup products $\phi \cup \psi$
	with $\phi \in H^1_{\stable}$ and $\psi \in H^1_{!,\stable}$.
	We have a linear map
	\[
		\Phi : V \to  \C,
		\qquad
		\Phi(\Sigma) = \langle \Sigma,\omega^{d-1}\rangle.
	\]
	The map $\Phi$ is a morphism of $\widehat{\G(\Q)}$ representations because $\omega$ is $\widehat{\G(\Q)}$-invariant.
	Using our assumption on $\Gamma'$, \autoref{lem:existence} shows that $\Phi$ is surjective.
	Therefore $V$ has a 1-dimensional trivial quotient.
	Since $H^2_{!,\stable}$ is semi-simple, $V$ must have
	a 1-dimensional trivial subrepresentation.
	By \autoref{lem:sigma stable},
	$(H^2_{!,\stable})^{\widehat{\G(\Q)}}$ is spanned by $\omega$.
	Therefore $\omega \in V$.
	In other words, there exist elements $\phi_i\in H^1_{\stable}$ and $\psi_i \in H^1_{!,\stable}$ such that
	\begin{equation}
		\label{eq:cup prod}
		\omega
		=
		\phi_1 \cup \psi_1 + \cdots + \phi_r \cup \psi_r.
	\end{equation}
	Choose an arithmetic subgroup $\Gamma''$ of sufficiently high level, so that all of the elements 	$\phi_i$ and $\psi_i$ are images in the direct limit of elements of $H^1(X'')$,
	where $X''= \Gamma'' \backslash \cH$.
	Then \eqref{eq:cup prod} holds in $H^2(X'')$.
	
	The group extension $\widetilde\SU(d,1)$ of $\SU(d,1)$ is represented by the cocycle $\sigma_\Z \in H^2_{\meas}(\SU(d,1),\Z)$.
	We shall write $\sigma_\C$ for the image of $\sigma_\Z$ in
	$H^2_{\cts}(\SU(d,1),\C)$.
	Recall that $\sigma_\C$ is a multiple of $\omega$.
	By (\ref{eq:cup prod})
	the restriction of $\sigma_\C$ to $\Gamma''$ is a sum
	of cup products of elements of $H^1(\Gamma'',\C)$.
	This means that the restriction of $\sigma_\Z$ to $\Gamma''$ satisfies the
	hypothesis of \autoref{lem:criterion}.
	By \autoref{lem:criterion}, the groups $\tilde \Gamma''$
	and $\tilde{\Gamma''}^{(n)}$ are all residually finite.
	Since $\tilde \Gamma \cap \tilde{\Gamma''}$ has finite index in $\tilde\Gamma$
	and $\tilde \Gamma \cap \tilde{\Gamma''}^{(n)}$ has finite index in $\tilde\Gamma^{(n)}$,
	it follows that $\tilde \Gamma$
	and $\tilde\Gamma^{(n)}$ are also residually finite.
\end{proof}

\section{Non-vanishing of the first inner cohomology}

In this section, we give an example of an arithmetic subgroup $\Gamma$
of $\SU(2,1)$ with cusps,
for which $H^1_!(X_\Gamma,\C) \ne 0$.
This demonstrates that the hypothesis of \autoref{thm:main} is satisfied in at least one case.
The construction which we describe here was suggested to the author by Matthew Stover,
and is a modification of his construction of the towers $C_j$ in section 5 of \cite{stover}.

We begin with the Deligne--Mostow group $\Gamma_\mu$, where
$\mu = (\frac{2}{6},\frac{2}{6},\frac{3}{6},\frac{4}{6},\frac{1}{6})$ (see \cite{deligne mostow}).
The group $\Gamma_\mu$ is an arithmetic subgroup od $\SU(2,1)$ with cusps (see the table on page 86 of \cite{deligne mostow}).
It is known (see section 5 of \cite{stover}) that there is a  finite index subgroup $\Gamma'\subseteq\Gamma_\mu$ for which there exists a surjective homomorphism
\[
	f : \Gamma' \to \Sigma,
\]
where $\Sigma$ is a hyperbolic surface group, i.e. the fundamental group of a compact Riemann surface of genus at least $2$.
The next theorem shows that whenever such a homomorphism $f$ exists, the hypothesis of \autoref{thm:main} is true.

\begin{theorem}
	\label{thm:surface}
	Let $\Gamma$ be an arithmetic subgroup of $\SU(d,1)$ with cusps.
	Assume that there exists a surjective homomorphism $\Gamma \to \Sigma$, where $\Sigma$ is a hyperbolic surface group.
	Then there exists an arithmetic subgroup $\Gamma' \subset \Gamma$, such that $H^1_!(X_{\Gamma'},\C) \ne 0$.
\end{theorem}

In the rest of this section, we shall prove \autoref{thm:surface} in a series of lemmata.
We shall assume from now on that
$\Gamma$ is an arithmetic subgroup of $\SU(d,1)$ with cusps, and that we have a surjective homomorphism $f: \Gamma \to \Sigma$, where $\Sigma$ is a hyperbolic surface group.
Replacing $\Gamma$ and $\Sigma$ by finite index subgroups if necessary, we shall assume that $\Gamma$ is torsion-free.

For a technical reason (in the proof of \autoref{lem:rho} below)
it will be more convenient in this section to work with
cohomology with coefficients in $\R$.

Recall that the cusps of $\Gamma$ are the elements $[v] \in \Proj^d(k)$ such that $(v,v)=0$, where $(-,-)$ is the Hermitian form of signature $(d,1)$.
In what follows, we shall abuse notation slightly by writing $v$ for a cusp, rather than $[v]$.
Let $v$ be a cusp of $\Gamma$, and let $\Gamma_v$ be the stabilizer of $v$ in $\Gamma$.
Since $\Gamma$ is torsion-free, the subgroup $\Gamma_v$ is
contained in the Heisenberg group $N_v$, and is a cocompact lattice in $N_v$.
In particular, the restriction map gives an isomorphism
$H^1_{\cts}(N_v,\R) \cong H^1(\Gamma_v,\R)$.

Given any subgroup $\Sigma'$ of finite index in $\Sigma$, there is
a linear map $R_{\Sigma',v} : H^1(\Sigma',\R) \to H^1_{\cts}(N_v,\R)$, defined as the following composition:
\[
	R_{\Sigma',v} : H^1(\Sigma' ,\R) \stackrel{f^*}{\to} H^1(\Gamma',\R) \stackrel{\Rest}\to H^1(\Gamma'_v,\R) \cong H^1_{\cts}(N_v,\R).
\]
Where $\Gamma' = f^{-1}(\Sigma')$.
If $\Sigma'' \subseteq \Sigma'$ is a subgroup of finite index, then
we have $R_{\Sigma',v}=R_{\Sigma'',v} \circ \Rest$, where $\Rest : H^1(\Sigma',\R) \to H^1(\Sigma'',\R)$ is the restriction map, so we actually have a map
\[
	R_v : \lim_{\to} H^1(\Sigma',\R) \to H^1_{\cts} (N_v,\R),
\]
where the direct limit is taken over all subgroups $\Sigma'$ of finite index in $\Sigma$.
Since the restriction maps $H^1(\Sigma',\R) \to H^1(\Sigma'',\R)$ are injective, the direct limit above may be regarded as a union
of an increasing sequence of finite dimensional vector spaces.

We shall call $v$ an \emph{essential cusp} is the map $R_v$ is non-zero.
This is equivalent to saying that
there exists a subgroup $\Sigma'$ of finite index in $\Sigma$,
such that the map $R_{\Sigma',v}$ is non-zero.
Note that if $R_{\Sigma',v}$ is non-zero then $R_{\Sigma'',v}$ is non-zero for each subgroup $\Sigma''\subseteq \Sigma'$.

\begin{lemma}
	\label{lem:conj Rv}
	Let $\Sigma'$ be a normal subgroup of finite index in $\Sigma$.
	If $R_{\Sigma',v}$ is non-zero, then for all $g \in \Gamma$
	the map $R_{\Sigma',gv}$ is non-zero.
\end{lemma}

\begin{proof}
	Assume that $R_{\Sigma',v}$ is non-zero.
	Choose an element $\phi \in H^1(\Sigma',\C)$ such that
	$R_{\Sigma',v}(\phi) \ne 0$.
	In other words, $\phi : \Sigma' \to \C$ is a homomorphism and the
	composition
	\[
		\Gamma'_v \hookrightarrow \Gamma' \stackrel{f}\to \Sigma' \stackrel{\phi}\to \C
	\]
	is non-zero.
	Choose an element $n \in \Gamma'_v$ whose image in $\C$ is non-zero.
	
	Define $\psi \in H^1(\Sigma',\C)$ by $\psi(\sigma) = \phi(f(g)^{-1} \sigma f(g))$.
	The element $n' = gng^{-1}$ is in $\Gamma'_{gv}$, and we have
	\[
		\psi( f ( n' )) = \phi(f(n)) \ne 0.
	\]
	Therefore $R_{\Sigma',gv}(\psi) \ne 0$.
\end{proof}

\begin{lemma}
	\label{lem:conj Rv 2}
	Let $\Sigma'$ be a subgroup of finite index in $\Sigma$
	and let $\Gamma'$ be the pre-image of $\Sigma'$ in $\Gamma$.
	For all $g\in\Gamma'$ and all cusps $v$ we have $\ker(R_{\Sigma',gv}) = \ker(R_{\Sigma',v})$.
\end{lemma}

\begin{proof}
	Let $\phi \in \ker (R_{\Sigma',v})$.
	This means that $\phi : \Sigma' \to \C$ is a homomorphism
	and $\phi(f(n))=0$ for all $n \in  \Gamma'_v$.
	If $n' \in \Gamma'_{gv}$ then we have $n' = gng^{-1}$ for some $n \in \Gamma'_v$.
	This implies
	\[
		\phi(f(n')) = \phi(f(g)) + \phi(f(n)) - \phi(f(g)) = 0.
	\]
	Hence $\phi \in \ker R_{\Sigma',gv}$.
	The converse is proved in the same way, replacing $g$ by $g^{-1}$.
\end{proof}

\begin{lemma}
	\label{lem:conj essential}
	If $v$ and $w$ are cusps and $w=gv$ for some $g \in \Gamma$
	then $v$ is essential if and only if $w$ is essential.
\end{lemma}

\begin{proof}
	This follows immediately from \autoref{lem:conj Rv}.
\end{proof}

\begin{lemma}
	\label{lem:sigma'}
	There exists a normal subgroup $\Sigma_0$ of finite index in $\Sigma$, such that for each essential cusp $v$, the map $R_{\Sigma_0,v}$ is non-zero.
\end{lemma}

\begin{proof}
	Let $v_1 , \ldots, v_r$ be a set of representatives for
	the $\Gamma$-orbits of the essential cusps.
	For each $i$, we may choose a normal subgroup $\Sigma_i$ of $\Sigma$, such that
	$R_{\Sigma_i,v_i}$ is non-zero.
	We shall prove the lemma with $\Sigma_0 = \Sigma_1 \cap \cdots \cap \Sigma_r$.
	Suppose $w$ is any essential cusp.
	By \autoref{lem:conj essential} we have $w=g v_i$ for some $g\in\Gamma$.
	By \autoref{lem:conj Rv} the map $R_{\Sigma_i,w}$ is non-zero.
	Since $\Sigma_0 \subseteq \Sigma_i$, it follows that $R_{\Sigma_0,w}$ is non-zero.
\end{proof}

\begin{lemma}
	\label{lem:rho}
	Let $\Sigma_0$ be chosen as in \autoref{lem:sigma'}
	and let $\Gamma_0$ be the pre-image of $\Sigma_0$ in $\Gamma$.
	There exists a surjective homomorphism $\rho : \Sigma_0 \to \Z$,
	such that for every essential cusp $v$, the composition
	$\rho \circ f : \Gamma_0 \to \Z$ is non-zero on $\Gamma_{0,v}$.
\end{lemma}

\begin{proof}
	Note that $\rho \circ f$ is non-zero on $\Gamma_{0,v}$ if and only if $R_{\Sigma_0,v}(\rho) \ne 0$.
	We must therefore show that there is an element of $H^1(\Sigma_0,\Z)$ which is not in the kernel of $R_{\Sigma_0,v}$ for any essential cusp $v$.

	Let $v_1 , \ldots, v_s$ be a set of representatives for
	the $\Gamma_0$-orbits of essential cusps.
	We have chosen $\Sigma_0$ so that each of the maps $R_{\Sigma_0,v_i}$ is non-zero.
	Hence $\ker R_{\Sigma_0,v_i}$ is a proper subspace of $H^1(\Sigma_0,\R)$.
	In particular, the union of the kernels of the $R_{\Sigma_0,v_i}$ is not the whole vector space $H^1(\Sigma_0,\R)$, and there is
	an open cone in $H^1(\Sigma_0,\R)$ which does not intersect any of these kernels.
	Choose a non-zero element $\rho \in H^1(\Sigma_0,\Z)$ in this open cone,
	so $\rho \not\in \ker R_{\Sigma_0,v_i}$ for all $i$.
	It follows from \autoref{lem:conj Rv 2}
	that $\rho \not\in \ker R_{\Sigma_0,v}$ for all essential cusps $v$.
	Dividing $\rho$ by a constant if necessary, we may assume that
	$\rho:\Sigma_0 \to \Z$ is surjective.
\end{proof}

Now let $\Sigma_0$ be chosen as in \autoref{lem:sigma'}
and let $\rho: \Sigma_0 \to \Z$ be a homomorphism chosen as
in \autoref{lem:rho}.
We now define a sequence of arithmetic groups $\Gamma_n$ as follows:
\[
	\Gamma_n = f^{-1}(\Sigma_n),
	\qquad\text{where }
	\Sigma_n = \{ \sigma \in \Sigma_0 : \rho(\sigma) \equiv 0 \bmod n \}.
\]

\begin{lemma}
	\label{lem:essential bound}
	The number of $\Gamma_n$-orbits of essential cusps is bounded
	independently of $n$.
\end{lemma}

\begin{proof}
	Choose any essential cusp $v$, and let $S_n$ be the set of
	$\Gamma_n$-orbits of cusps which are in the same $\Gamma_0$-orbit as $v$.
	It is sufficient to show that the cardinality of each $S_n$ is bounded independently of $n$.
	By the orbit-stabilizer theorem there is a bijection
	between $S_n$ and the double coset set:
	\[
		 S_n \simeq \Gamma_n \backslash \Gamma_0 / \Gamma_{0,v}.
	\]
	Using the homomorphism $\rho \circ f$, we may identify
	$\Gamma_n \backslash \Gamma_0$ with $\Z/n\Z$.
	Therefore there is a bijection between $S_n$ and the group $\Z / (n\Z + \rho(f(\Gamma_{0,v}))$.
	In particular we have
	\[
		| S_n | \le | \Z / \rho (f(\Gamma_{0,v})) |.
	\]
	The homomorphism $\rho$ is chosen so that $\rho (f(\Gamma_{0,v})) \ne 0$, so we have a bound on the cardinality of $S_n$ which does not depend on $n$.
\end{proof}

\begin{lemma}
	\label{lem:rank bound}
	The rank of the composition $H^1(\Sigma_n,\R) \stackrel{f^*}\to H^1(\Gamma_n,\R) \to H^1(\partial X_{\Gamma_n},\R)$ is bounded independently of $n$.
\end{lemma}

\begin{proof}
	There is a decomposition
	\[
		H^1(\partial X_{\Gamma_n},\R)
		=
		\bigoplus_v H^1(\Gamma_{n,v},\R),
	\]
	where $v$ ranges over the $\Gamma_n$-orbits of cusps.
	However, if $v$ is not an essential cusp, then the
	map $H^1(\Sigma_n,\C) \to H^1(\Gamma_{n,v},\C)$ is zero.
	Therefore, the image of $H^1(\partial X_{\Gamma_n},\C)$ is contained in the direct sum of the spaces $H^1(\Gamma_{n,v},\C)$, where $v$ ranges over the $\Gamma_n$-orbits of the essential cusps.
	The result now follows from \autoref{lem:essential bound}.
\end{proof}

\begin{lemma}
	\label{lem:last}
	For $n$ sufficiently large, the inner cohomology $H^1_!(X_{\Gamma_n},\R)$ is non-zero.
\end{lemma}

\begin{proof}
	Since $\Sigma_0$ is a hyperbolic surface group,
	the dimension of $H^1(\Sigma_n,\R)$ tends to infinity as $n \to \infty$.
	In view of \autoref{lem:rank bound}, for large enough $n$, the map $H^1(	\Sigma_n,\R) \to H^1(\partial X_{\Gamma_n},\R)$ is not injective.
	Choose a non-zero element $\phi \in \ker (H^1(	\Sigma_n,\R) \to H^1(\partial X_{\Gamma_n},\R) )$.
	Then $f^*(\phi)$ is a non-zero element of $H^1_!(\Gamma_n,\R)$.
\end{proof}

\autoref{lem:last} concludes the proof of \autoref{thm:surface}.
By the discussion above, this shows that the preimage of 
$\Gamma_{(\frac{2}{6},\frac{2}{6},\frac{3}{6},\frac{4}{6},\frac{1}{6})}$ in each every connected cover of $\SU(2,1)$ is residually finite.

\end{document}